\documentclass[preprint,12pt]{elsarticle}

\usepackage{enumitem}
\usepackage{algorithm}
\usepackage{algorithmic}
\usepackage{subfig,graphicx,amsmath,latexsym,amssymb,amsthm}
\usepackage[colorlinks=true]{hyperref}
\usepackage{subfig,graphicx,amsmath,latexsym,amssymb,amsthm}
\usepackage[colorlinks=true]{hyperref}
\usepackage{float}
\usepackage{booktabs}
\usepackage{bm}
\usepackage[labelfont=bf]{caption} 

\usepackage{xcolor}

\journal{}

\begin{document}
	
	\begin{frontmatter}
		
		
		
		\title{Deep Learning-Enhanced Calibration of the Heston Model: A Unified Framework} 
		
		\author[guilan]{Arman Zadgar}
		\author[alzahra]{Somayeh Fallah}
		\author[guilan]{Farshid Mehrdoust\corref{cor1}}
		\cortext[cor1]{Corresponding author}
		\ead{fmehrdoust@guilan.ac.ir} 
		\author[Almeria]{Juan E. Trinidad Segovia}
		
		\address[guilan]{Department of Applied Mathematics, Faculty of Mathematical Sciences,
			University of Guilan, P. O. Box: 41938-1914, Rasht, Iran}
		
		\address[alzahra]{Department of Mathematics, Faculty of Mathematical Sciences, Alzahra University, Tehran, Iran}
		
		\address[Almeria]{Department of Economics and Business, Universidad de Almeria, Spain}
		
		\begin{abstract}
			The Heston stochastic volatility model is a widely used tool in financial mathematics for pricing European options. However, its calibration remains computationally intensive and sensitive to local minima due to the model's nonlinear structure and high-dimensional parameter space. This paper introduces a hybrid deep learning-based framework that enhances both the computational efficiency and the accuracy of the calibration procedure. The proposed approach integrates two supervised feedforward neural networks: the Price Approximator Network (PAN), which approximates the option price surface based on strike and moneyness inputs, and the Calibration Correction Network (CCN), which refines the Heston model's output by correcting systematic pricing errors. Experimental results on real S\&P 500 option data demonstrate that the deep learning approach outperforms traditional calibration techniques across multiple error metrics, achieving faster convergence and superior generalization in both in-sample and out-of-sample settings. This framework offers a practical and robust solution for real-time financial model calibration.
		\end{abstract}
		
		
		
		
	\end{frontmatter}
	
		
		
		\section{Introduction}
		The Heston stochastic volatility model, first introduced by \cite{hes}, has established itself as a fundamental tool in financial mathematics, offering a more realistic framework for pricing European options than earlier models such as that of \cite{bls}. By incorporating stochastic volatility--unlike the constant volatility assumption in the Black--Scholes framework--the Heston model is capable of capturing important empirical features observed in financial markets, including volatility clustering and the volatility smile. Its semi-closed-form solution also supports computational efficiency, which is crucial for real-time trading and risk management applications \cite{gath, duff}. However, the calibration of the Heston model remains a non-trivial task due to the complex and nonlinear relationships among its parameters, particularly under conditions of heightened market uncertainty. This complexity highlights the growing need for advanced optimization algorithms and machine learning techniques to achieve accurate and robust calibration \cite{Broadie, and, Gneiting}.
		
		In essence, model calibration can be framed as an inverse problem, wherein the objective is to infer the underlying model parameters from observed outputs--such as market option prices. A rich body of research has explored a variety of techniques to address this inverse problem, including adjoint optimization, Bayesian inference, and sparsity-regularized methods.
		
		In financial practice, particularly in the context of derivative pricing and risk management, asset model calibration involves recovering the parameters of stochastic differential equations (SDEs) that govern asset dynamics using real market data. For equity and options markets, this means identifying the appropriate model parameters so that the mathematical model reproduces the observed prices of actively traded options with high fidelity. Once calibrated, these models become essential tools for pricing exotic over-the-counter instruments and for managing hedging and portfolio risk in a dynamic environment.
		
		Given the frequency with which models must be recalibrated in professional settings--often several times per day--the calibration process must be not only accurate, but also computationally efficient and stable. The demand for real-time decision-making imposes strict performance requirements, particularly when dealing with high-dimensional models or complex financial products, making traditional calibration methods increasingly impractical in modern market conditions.
		
		In recent years, deep learning has emerged as a transformative tool in finance, offering sophisticated methods to analyze and predict the complex behaviors of financial markets. Leveraging advanced algorithms and computational resources, deep learning models process vast amounts of data with unprecedented efficiency, with neural networks excelling at capturing nonlinear relationships often missed by traditional models \cite{heat, good, gu}. Such capabilities are indispensable for activities such as stock price forecasting, portfolio optimization, and arbitrage identification. Deep learning also complements classical models like Black--Scholes and Heston by approximating complex option price surfaces more accurately and addressing residual errors from simplifying assumptions, significantly reducing computational burdens in methods like Monte Carlo simulations \cite{hut, siri, ruf}. The integration of machine learning with traditional models represents a paradigm shift, enabling robust, adaptable pricing frameworks that enhance risk management and derivative pricing under evolving market conditions, fostering financial innovation and deeper insights into market behavior \cite{bueh}.
		
		Motivated by these developments, this study proposes a novel approach that integrates model calibration with deep learning, viewing the problem through the lens of optimization. By harnessing the approximation power and generalization capabilities of deep neural networks, we construct a data-driven calibration framework that captures the complex, nonlinear relationships among model parameters, option features, and observed market prices. This paradigm shift offers both computational efficiency and robustness--two critical requirements for modern financial modeling.
		
		Specifically, we present a hybrid calibration method that augments the Heston stochastic volatility model using two regression-based neural networks trained in a supervised learning setting. The first, the PAN, approximates the relationship between strike prices and option values, providing a smooth surface for in-sample and out-of-sample pricing \cite{ford, Chri}. The second, the CCN, corrects the residual discrepancies between the Heston model's output and observed market data \cite{Li, Chen}. By cascading these networks, we create a two-phase system that bridges traditional stochastic modeling with modern machine learning \cite{Hor, Beck}. This unified framework achieves high pricing accuracy while maintaining the interpretability and tractability of the Heston model, positioning it as a robust tool for real-time financial applications.
		
		The remainder of this paper is structured as follows. Section~\ref{sec:calibration} outlines the theoretical formulation of the Heston stochastic volatility model and presents the associated calibration framework. Section~\ref{sec:nn_model} introduces the use of neural networks in model calibration, emphasizing their role as surrogate pricing functions. Section~\ref{sec:ffnn} describes the architecture of feedforward neural networks, while Section~\ref{sec:framework} discusses the computational environment used for their implementation. Section~\ref{sec:training} details the training procedures and optimization techniques employed to ensure robust and efficient learning. Section~\ref{sec:implementation} illustrates the integration of deep learning into the Heston model calibration process through the design of the Price Approximator Network (PAN) and the Calibration Correction Network (CCN). Section~\ref{results} presents empirical results and evaluates the performance of the proposed method in comparison with traditional calibration approaches. Finally, Section~\ref{conclusion} concludes the paper with a summary of contributions and suggestions for future research directions.

		\section{Financial model calibration} \label{sec:calibration}
		In the realm of financial mathematics, understanding and modeling market behaviors accurately is paramount. The marriage of traditional stochastic models with cutting-edge deep learning techniques offers a robust framework for tackling the complexities of modern financial markets. This section delves into the theoretical underpinnings that guide our hybrid approach, starting with a detailed discussion of Heston's Stochastic Volatility Model. Integrating these core concepts, our objective is to facilitate a more nuanced and precise approach to the calibration of option pricing models.
		
		\subsection{Heston's Stochastic Volatility Model}
		The Heston model is a a widely used stochastic volatility model in financial mathematics that enhances the Black--Scholes model by accounting for the stochastic nature of volatility. This model is particularly effective in capturing empirical phenomena like the volatility smile and volatility clustering, which are often observed in real world in financial markets.
		
		In the Heston model, the contribution refers to the treatment of a financial market characterized by a finite maturity horizon \( T \) and a risky asset $S=\{ S(t), \ 0 \leq t \leq T \}$ whose stochastic price dynamic is defined, over the probability space $ (\Omega, \mathcal{F}, \mathbb{F} = \{ \mathcal{F}_t, \ t \geq 0 \}, \mathbb{P})$ by the following system of stochastic differential equation:
		\begin{equation*}
			dS(t) = r S(t) dt + \sqrt{v(t)} S(t) dW_1(t),
		\end{equation*}
		where \( r \) is a constant risk--free interest rate of the asset, \( v(t) \) is the instantaneous variance, encapsulating the asset's volatility, and \( W_1(t) \) is a Brownian motion process, driving the asset price.
		
		The volatility \( v(t) \) is modeled as a mean-reverting process, specifically a Cox-Ingersoll-Ross (CIR) process \cite{Cox, Stein}:
		\begin{equation*}
			dv(t) = \kappa (\theta - v(t)) dt + \sigma \sqrt{v(t)} dW_2(t),
		\end{equation*}
		where \( \kappa \) is the rate at which \( v(t) \) reverts to its long-term mean \( \theta \); \( \theta \) is the long-term mean level of variance; \( \sigma \) is the volatility of the volatility process; and \( W_2(t) \) is another Brownian motion process, correlated with \( W_1(t) \) through a correlation coefficient \( -1 \leq \rho \leq 1 \). The correlation \( \rho \) between \( W_1(t) \) and \( W_2(t) \) is essential for capturing the leverage effect, which indicates that negative returns are typically associated with rising volatility, a phenomenon observable in actual financial markets.

		To price options under the Heston model, we derive the governing PDE by applying It\^o's lemma to the option price function \( C(S, v, t) \). According to It\^o's lemma, the differential \( dC \) of a sufficiently smooth function \( C(S, v, t) \), where \( S \) and \( v \) are stochastic processes, is expressed as follows:
		\begin{align*}
			dC &= \frac{\partial C}{\partial t} dt + \frac{\partial C}{\partial S} dS + \frac{\partial C}{\partial v} dv + \frac{1}{2} \frac{\partial^2 C}{\partial S^2} dS^2 + \frac{1}{2} \frac{\partial^2 C}{\partial v^2} dv^2 + \frac{\partial^2 C}{\partial S \partial v} dS dv.
		\end{align*}
		
		Using the SDEs for \( S(t) \) and \( v(t) \), we have:
		\begin{align*}
			dS &= rS dt + \sqrt{v} S dW_1, \\
			dS^2 &= (\sqrt{v} S)^2 dt = v S^2 dt, \\
			dv &= \kappa(\theta - v) dt + \sigma \sqrt{v} dW_2, \\
			dv^2 &= (\sigma \sqrt{v})^2 dt = \sigma^2 v dt, \\
			dS dv &= \rho \sigma v S dt.
		\end{align*}
		
		Substituting these into It\^o's lemma, we get:
		\begin{align*}
			dC &= \frac{\partial C}{\partial t} dt + \frac{\partial C}{\partial S} (rS dt + \sqrt{v} S dW_1) + \frac{\partial C}{\partial v} (\kappa(\theta - v) dt + \sigma \sqrt{v} dW_2) \\
			&\quad + \frac{\partial^2 C}{\partial S^2} \frac{1}{2} v S^2 dt + \frac{\partial^2 C}{\partial v^2} \frac{1}{2} \sigma^2 v dt + \frac{\partial^2 C}{\partial S \partial v} \rho \sigma v S dt.
		\end{align*}
		
		Rearranging terms, we obtain the PDE:
		\begin{align*}
			\frac{\partial C}{\partial t} + rS \frac{\partial C}{\partial S} + \kappa(\theta - v) \frac{\partial C}{\partial v} + \frac{1}{2} v S^2 \frac{\partial^2 C}{\partial S^2} + \rho \sigma v S \frac{\partial^2 C}{\partial S \partial v} + \frac{1}{2} \sigma^2 v \frac{\partial^2 C}{\partial v^2} - rC = 0,
		\end{align*}
			with the terminal boundary condition for a European call option:
			\begin{equation*}
				C(S_T, v, T) = \max(S_T - K, 0),
			\end{equation*}
			where $T$ denotes the maturity date, $S_T$ is the underlying asset price at maturity, and $K$ represents the strike price. 
			Note that the payoff depends solely on the terminal asset price and is independent of the variance state variable.
		
		One of the powerful aspects of the Heston model is its ability to provide a closed form solution for the price of a European call option using the characteristic function. The characteristic function \( f(\phi; x_t, v_t) \) is given by:
		\begin{align*}
			f(\phi; x_t, v_t) = \exp\left( i\phi x_t + C(\tau, \phi) + D(\tau, \phi)v_t \right),
		\end{align*}
		where \( x_t = \log(S_t) \) and \( \tau = T - t \).
		
		To derive the characteristic function, we start from the general solution of the Heston model and the properties of the characteristic function. First, we define the logarithm of the asset price as \( x_t = \log(S_t) \) and formulate the characteristic function as:
		\begin{align*}
			f(\phi; x_t, v_t) = \mathbb{E}\left[e^{i\phi x_T} \mid x_t, v_t\right].
		\end{align*}
		
		Using the Feynman--Kac theorem, the characteristic function can be expressed as:
		\begin{align*}
			f(\phi; x_t, v_t) = \exp\left( i\phi x_t + C(\tau, \phi) + D(\tau, \phi)v_t \right).
		\end{align*}
		
		The functions \( C(\tau, \phi) \) and \( D(\tau, \phi) \) are solutions to the following system of Riccati differential equations:
		\begin{align*}
			\frac{\partial C}{\partial \tau} &= \kappa\theta D - \frac{\sigma^2}{2} D^2 + \frac{1}{2} (i\phi - \phi^2), \\
			\frac{\partial D}{\partial \tau} &= \kappa D - \rho\sigma i\phi - \frac{\sigma^2}{2} D^2.
		\end{align*}
		
		By integrating these Riccati equations, we obtain:
		\begin{align*}
			C(\tau, \phi) &= \frac{\kappa \theta}{\sigma^2} \left[ (\kappa - \rho\sigma i\phi)\tau - 2\log\left( \frac{1 - g e^{-d \tau}}{1 - g} \right) \right], \\
			D(\tau, \phi) &= \frac{\kappa - \rho\sigma i\phi + d}{\sigma^2} \left[ \frac{1 - e^{-d \tau}}{1 - g e^{-d \tau}} \right],
		\end{align*}
		with auxiliary functions:
		\begin{align*}
			d &= \sqrt{(\rho\sigma i\phi - \kappa)^2 - \sigma^2 (i\phi + \phi^2)}, \\
			g &= \frac{\kappa - \rho\sigma i\phi + d}{\kappa - \rho\sigma i\phi - d}.
		\end{align*}
		
		Using these definitions, the price of a European call option \( C(S_0, K, T) \) is obtained by the inverse Fourier transform of the characteristic function:
		\begin{align*}
			C(S_0, K, T) = S_0 P_1 - K e^{-rT} P_2,
		\end{align*}
		where \( P_1 \) and \( P_2 \) represent the risk-neutral probabilities of the option finishing in the money:
		\begin{align*}
			P_1 &= \frac{1}{2} + \frac{1}{\pi} \int_0^\infty \text{Re} \left[ \frac{e^{-i\phi \log K} f(\phi - i)}{i\phi f(-i)} \right] d\phi, \\
			P_2 &= \frac{1}{2} + \frac{1}{\pi} \int_0^\infty \text{Re} \left[ \frac{e^{-i\phi \log K} f(\phi)}{i\phi} \right] d\phi.
		\end{align*}
		
		The Heston model's characteristic function provides a robust framework for deriving closed-form solutions for European call options. By leveraging the characteristic function, we efficiently calculate option prices and understand the impact of stochastic volatility on option valuation. This detailed explanation bridges theoretical models with practical applications, offering valuable insights for both researchers and practitioners in financial markets.

		\subsection{Calibration Framework}
		To calculate option prices, certain model parameters are needed that aren't directly observable from market data. The process of adjusting these parameters so that the model's predicted prices align with actual market prices is known as calibration. A major difficulty in this process is that market data alone often provides insufficient information to accurately infer the model parameters. In fact, different combinations of parameters can result in model prices that appear consistent with the market.
		
		In reality, achieving a perfect match between model and market prices is not feasible. As a result, calibration becomes an optimization problem. The aim is to minimize an objective function--also known as an error function--that quantifies the difference between model-generated prices and market prices for a given parameter set.
		
		There are various ways to define this objective function in the literature. In this work, the Weighted Root Mean Square Error (Weighted RMSE) is used to measure the difference between model and market prices, based on a parameter set \( \mathcal{H} \).
		
		Let \( \eta \in \mathcal{H} \) denote the set of model parameters, \( m = \ln(S_t / K) \) represent the log-moneyness (where \( K \) is the strike price and \( S_t \) is the current spot price), and \( T \) denote the option's time to maturity. Let \( \hat{V} \) and \( V^{\text{Mkt}} \) represent the model-generated and market-observed option prices, respectively. The weighted RMSE is then defined as:
		\begin{equation*} 
			J(\eta) = \sqrt{ \sum_i \sum_j \omega_{i,j} \left( \hat{V}(\eta, T_i, m_j) - V^{\text{Mkt}}(T_i, m_j) \right)^2 }.
		\end{equation*}
		Here, \( \omega_{i,j} \) are the weights assigned to each option price, which depend on the option's log-moneyness and maturity. These weights are chosen to reflect both the relative importance of each option and the reliability of its observed market price.
		
		To address the optimization challenge involved in the calibration process, a variety of methods have been proposed. Two notable techniques include:
		
		\begin{itemize}
			\item \textbf{Nelder-Mead Method:}
			The Nelder--Mead algorithm, introduced by John Nelder and Roger Mead in 1965, is a heuristic method for solving nonlinear optimization problems. It aims to minimize a continuous function in a multi-dimensional space without requiring derivative information. Also known as the \textit{downhill simplex method}, it operates on a geometric structure called a simplex---a polytope with \( N + 1 \) vertices in an \( N \)-dimensional space. The algorithm begins with an initial simplex and then iteratively transforms it by stretching, shrinking, and moving it toward the region where the function reaches a local minimum.
			
			\item \textbf{Differential Evolution:}
			Differential Evolution (DE) is another powerful optimization method that does not rely on gradient information and doesn't require prior parameter initialization. It's especially useful for finding global minima in non-convex objective functions. In this method, as described by \cite{Liu}, a population of candidate solutions is generated. For each individual \( \eta_i \), a mutant \( \eta_i' \) is created using:
			\begin{equation*}
				\eta_i' = \eta_a + F \cdot (\eta_b - \eta_c),
			\end{equation*}
			where \( a \neq i \), and \( a, b, c \) are randomly selected indices from the population. The parameter \( F \in [0, \infty) \) is known as the differential weight and controls the mutation step size. Depending on the strategy (e.g., \texttt{rand1bin} or \texttt{best1bin}), \( \eta_a \) may be chosen randomly or as the best solution in the previous generation. A crossover process filters the candidates based on a probability \( Cr \).
			If the mutant \( \eta_i' \) yields a lower value of the objective function \( J \), it replaces the original candidate \( \eta_i \) in the population.
			
			This process is repeated until convergence or a stopping condition is met. The effectiveness of the algorithm heavily depends on parameter tuning. Larger mutation factors and population sizes can improve the chances of locating a global minimum. Additionally, the convergence tolerance parameter, which measures population diversity, helps determine when to terminate the optimization.
			
		\end{itemize}
		
		The speed and effectiveness of the calibration process are heavily influenced by the computational cost of pricing vanilla options, as these evaluations are repeatedly required during the optimization of the objective function. Traditional approaches, while theoretically sound, often become impractical due to the high dimensionality and complexity of the Heston model. Inspired by the work of \cite{Liu} and \cite{Hor}, this study proposes the use of Artificial Neural Networks (ANNs) as surrogate models to approximate the pricing function. This substitution enables significant acceleration of the calibration routine, offering a more efficient solution without sacrificing accuracy. By doing so, we transition from computationally expensive numerical methods to a learning-based framework capable of handling large-scale data and complex parameter interactions with enhanced robustness and speed.

		\section{Neural Networks in Model Calibration} \label{sec:nn_model}
		Modern financial model calibration faces dual challenges: accurately capturing market dynamics while maintaining computational tractability. Traditional numerical methods often struggle to balance these requirements, particularly for high-dimensional models with non-linear dependencies. Neural networks emerge as a transformative paradigm in this context, leveraging their universal approximation capabilities to simultaneously address accuracy and efficiency constraints. By learning complex mappings between model parameters and financial instruments' prices, neural networks enable robust calibration frameworks that adapt to diverse market regimes. This section examines their foundational role, beginning with the core mechanism of function approximation.
		
		\subsection{Neural Network Function Approximation}
		Suppose we aim to estimate a true pricing function \( f(x) \) using an approximate function \( \hat{f}(x) \). Neural networks, denoted as \( F(x, W) \), where \( x \) represents the input and \( W \) the set of weights, can effectively approximate such functions. When trained on a data set of input-output pairs \( (x, y = f(x)) \), the network learns optimal weights \( \hat{W} \), resulting in the approximation:
		\begin{equation*}
			F(x, \hat{W}) = \hat{f}(x).
		\end{equation*}
		
		In recent years, artificial neural networks (ANNs) have been increasingly applied to model calibration tasks. A notable method proposed by A.~Hernandez treats calibration as an inverse mapping problem--from the market implied volatility surface to the corresponding model parameters. Hernandez applied this technique to calibrate the Hull--White model. This approach is particularly appealing, as ANNs can directly infer model parameters from implied volatility surfaces, potentially bypassing iterative numerical optimization. 
		
		However, a major drawback is its reliance on historical implied volatility data for training. Since high-quality historical data is often limited, the model may overfit and perform poorly on new, unseen market conditions--especially in scenarios involving regime shifts. As a result, frequent retraining of the neural network is necessary, which is both computationally intensive and time-consuming.
		
		To address these challenges, a two-step alternative approach has been proposed. This method has been applied by \cite{Bayer} for calibrating the Rough Bergomi model, and by \cite{Liu} for the Heston and Bates models.
		
		\subsection{Understanding Calibration Slowness}
		This section explores the reasons behind the slow performance of the calibration process. Typically, thousands of pricing function evaluations are needed for options with various strikes and maturities to determine the best-fit parameters. However, the pricing function is not known in closed form and must be approximated using numerical techniques such as Monte Carlo simulations, Fourier transforms, or the resolution of partial differential equations through finite difference methods.
		
			As a result, the minimized objective function is formulated as
			\begin{equation*}
				J(\eta) = 
				\sum_{i=1}^{N_T} \sum_{j=1}^{N_m} 
				\omega_{i,j} 
				\left(
				\hat{V}(\eta, T_i, m_j) - V^{\text{Mkt}}(T_i, m_j)
				\right)^2,
			\end{equation*}
			where $T_i$ denotes the $i$-th maturity ($i = 1,\dots,N_T$), $m_j$ denotes the $j$-th log-moneyness level ($j = 1,\dots,N_m$), and $\omega_{i,j}$ represents the corresponding weighting coefficient. The term $\hat{V}(\eta, T_i, m_j)$ refers to the model-implied option price computed--by numerical or semi-analytical methods--for the specific maturity-moneyness pair $(T_i, m_j)$ under parameter vector $\eta$, while $V^{\text{Mkt}}(T_i, m_j)$ denotes the observed market price for the same contract.

		These numerical methods are computationally demanding and represent a primary factor contributing to the slow calibration process. Additionally, they may introduce numerical instability due to inherent approximation errors. 
		
		To address these issues, neural networks can be employed as a substitute for the numerical pricing step. This replacement transforms the optimization into a more deterministic and significantly faster process. Instead of recalculating option prices using expensive numerical simulations, the objective function uses outputs from a pre-trained artificial neural network (ANN), dramatically enhancing calibration efficiency.
		
		\subsection{Approximating the Pricing Function Using ANN}
		
		Let \( \hat{V} \) represent the model-implied option price, parameterized by \( \eta \) (the set of model parameters), the time to maturity \( T \), and the log-moneyness \( m \).
		
		To replicate this pricing function, we employ a supervised learning strategy where a neural network is trained to approximate the mapping from inputs \( (\eta, T, m) \) to the corresponding option prices \( \hat{V} \). The learning process consists of optimizing the network parameters \( W \) by minimizing a loss function over a synthetic data set of size \( N \). This is typically done by solving the following optimization problem based on the Weighted Root Mean Squared Error (RMSE):
		\begin{equation*}
			\hat{W} := \arg\min_{W \in \mathbb{R}^n} \sqrt{\sum_{i=1}^N \left( F(W, \eta_i, T_i, m_i) - \hat{V}(\eta_i, T_i, m_i) \right)^2}.
		\end{equation*}
		Here, \( F(W, \cdot) \) denotes the neural network function parameterized by weights \( W \). Optimization is commonly performed using first-order methods such as SGD or more sophisticated variants like the Adam algorithm. Implementation specifics are detailed in Section \ref{sec:training}.
		
		\subsection{Using the Trained Neural Network for Model Calibration}
		After successfully training the neural network and obtaining the optimal weights \( \hat{W} \), the network serves as a fast and accurate proxy for the original pricing function. This enables an efficient reformulation of the calibration problem: rather than repeatedly solving a computationally expensive pricing model, we use the surrogate network to estimate prices and calibrate the model parameters \( \eta \) by minimizing the squared deviation from observed market prices \( V^{\text{Mkt}} \):
		\begin{equation*}
			\hat{\eta} := \arg\min_{\eta \in \mathcal{H}} \sum_i \sum_j \omega_{i,j} \left( F(\hat{W}, \eta_i, T_i, m_j) - V^{\text{Mkt}}(T_i, m_j) \right)^2.
		\end{equation*}
		The training process requires the generation of a synthetic data set, capturing a representative range of model configurations and market conditions to ensure the neural network accurately learns the pricing manifold.
		
		\vspace{.3cm}
		
		The approach adopted in this research follows a robust and systematic two-phase framework designed to leverage the predictive capabilities of neural networks for model calibration:
		\begin{itemize}
			\item[1.] Construct a synthetic data set encompassing a comprehensive range of input parameters related to the financial model, option specifications, and prevailing market conditions.
			\item[2.] Calculate the corresponding option prices using the chosen pricing methodology, establishing the ground truth for supervised learning.
			\item[3.] Partition the data set into training and validation subsets to ensure reliable performance assessment.
			\item[4.] Train the ANN on the training set to approximate the complex pricing function, and evaluate its generalization accuracy on the test set.
			\item[5.] Utilize the trained neural network to carry out the calibration process, enabling efficient estimation of model parameters based on observed market data.
		\end{itemize}
		
		The proposed methodology employs a supervised deep learning regression model to approximate the pricing function, where the neural network learns to map input parameters to corresponding option prices for efficient model calibration.
		This methodology demonstrates strong stability and reliability. By training neural networks on synthetically generated data rather than relying on limited historical data sets, the model gains the ability to generalize effectively to unseen and future market scenarios. This anticipatory capability significantly enhances its practical relevance in dynamic financial environments.
		
		A key advantage of this approach is its ability to disentangle the overall calibration error into two components: the approximation error introduced by the neural network and the model's deviation from actual market data. Because the neural network requires training only once, the overall process becomes more robust and computationally efficient in the long run.
		
		Notably, in \cite{green}, it is observed that neural networks yield higher accuracy compared to traditional Monte Carlo methods. Furthermore, the use of backpropagation enables rapid and accurate evaluation of both the neural network output and its gradient with respect to the model parameters. This results in significantly faster option pricing and model calibration. As previously emphasized, this efficiency arises because the ANN transforms the inherently stochastic optimization problem into a deterministic one, streamlining the entire calibration process.

		\section{Feedforward Neural Network Architecture} \label{sec:ffnn}
		Deep learning, a subfield of machine learning, focuses on the training of multi-layer neural networks capable of capturing complex and nonlinear relationships within data. In this study, we employ feedforward neural networks (FFNNs) as the foundational architecture. However, due to the depth of the network and the scale of the data utilized, our approach is situated within the deep learning paradigm. The distinction between traditional neural networks and deep learning primarily lies in the depth of the model, the volume of data required, and the computational resources involved--such as the use of GPUs for accelerated training.
		
		Although FFNNs represent the foundational architecture in the realm of neural networks, they are recognized for their remarkable capacity to approximate complex, high-dimensional functions. In this study, the neural network employed is a standard FFNN, and this section outlines its structural and functional characteristics.
		
		Let us consider a neural network comprising \( L \) layers, where each layer is indexed by \( l \in \{1, \ldots, L\} \). Denote the input vector by \( x \in \mathbb{R}^n \) and the network's output by \( y \in \mathbb{R} \). For each layer \( l \), let \( W^{(l)} \in \mathbb{R}^{n_{l-1} \times n_l} \) and \( b^{(l)} \in \mathbb{R}^{n_l} \) represent the weight matrix and bias vector, respectively.
		
		The forward propagation through the network is defined recursively as follows:
		\begin{align*}
			z^{(0)} &= x, \\
			z^{(l)} &= \phi\left(z^{(l-1)}\right) W^{(l)} + b^{(l)} \quad \text{for } l = 1, \ldots, L, \\
			y &= z^{(L)},
		\end{align*}
		where \( z^{(l)} \in \mathbb{R}^{n_l} \) is the vector of neuron activations in layer \( l \), and \( \phi: \mathbb{R} \rightarrow \mathbb{R} \) is a nonlinear activation function applied element-wise. The nonlinearity of \( \phi \) is crucial; without it, multiple-layer networks could be reduced to a single linear transformation, negating the representational power of deeper architectures.
		
		The first and last layers correspond to the input (\( l = 0 \)) and output (\( l = L \)) layers, respectively. All layers in between, indexed by \( l \in \{1, \ldots, L-1\} \), are referred to as hidden layers, where the network learns to model intricate relationships within the data.
		
		To operationalize this architecture in a practical setting, it is necessary to adopt a computational framework that supports both flexibility in model design and efficiency in execution. In this work, we utilize one of the most advanced and widely used platforms for implementing deep learning models. The following subsection provides a detailed discussion of the selected framework and its relevance to our regression-based calibration strategy.

		\section{Framework for Neural Network Modeling} \label{sec:framework}
		Deep learning regression models can be implemented through a range of advanced computational frameworks, among which TensorFlow and PyTorch stand out as the most widely adopted and robust platforms for research and industrial deployment \cite{Plat, idmeh}. In the present study, PyTorch was selected for the development and training of our regression models due to its exceptional flexibility and widespread acceptance within the machine learning community. Known for its imperative, Pythonic programming interface, PyTorch provides an intuitive and accessible environment that facilitates rapid prototyping and experimentation, making it particularly appealing to researchers and developers.
		
		One of PyTorch's defining features is its use of dynamic computation graphs, which allow for real-time modification of network architecture during the training process. This capability affords enhanced control and significantly improves debugging workflows compared to static graph-based frameworks, thereby supporting more iterative and adaptive model development. In addition, PyTorch offers high computational efficiency through seamless integration with GPU acceleration, enabling the rapid execution of large-scale models and effectively addressing the intensive computational demands of deep neural networks.
		
		Beyond its core functionalities, PyTorch is supported by a comprehensive ecosystem of libraries and tools that streamline every stage of the machine learning pipeline--from data preprocessing and model visualization to deployment and inference. Its extensive collection of pre-built modules and utilities simplifies the construction of sophisticated network architectures, making it a preferred framework for both academic inquiry and real-world application.
		
		Taken together, these advantages make PyTorch an ideal platform for implementing our deep learning regression models. Its dynamic structure, computational performance, and developer-friendly environment enable us to achieve a high degree of accuracy, scalability, and efficiency in our modeling tasks.
		
		\vspace{0.3cm}
		Having defined the architecture and established the deep learning framework for model construction, we now turn to the training and optimization process. Effective training is essential to unlock the predictive power of neural networks and ensure they generalize well beyond the training data. The next section details the strategies used to initialize weights, propagate errors, define loss functions, and optimize model parameters using both classical and adaptive gradient-based methods.

		\section{Training Strategies and Optimization Algorithms} \label{sec:training}
		Having established the architectural foundation and implementation environment of the Feedforward Neural Network (FFNN), we now turn our attention to the training and optimization procedures that underpin its predictive capabilities. While architecture determines the expressive power of the model, effective training strategies are essential to realize this potential and ensure convergence to an accurate solution.
		
		This section presents the methodologies used to train the neural networks introduced in the previous sections. Topics include weight initialization techniques, the mechanics of forward and backward propagation, cost function formulation, and gradient-based optimization strategies. Special emphasis is placed on stochastic methods such as SGD and adaptive algorithms like Adam, which are instrumental in achieving robust and efficient model calibration in high-dimensional financial settings.
		
		These optimization tools not only accelerate the learning process but also enhance model generalization, making them integral to the success of our deep learning-based calibration framework for the Heston model.
		
		\subsection{Weight Initialization and Kaiming Method}
		The effectiveness of training a neural network is highly influenced by its initial weight configuration. Consequently, proper weight initialization is a crucial step before commencing the training process. While basic strategies often involve random initialization using uniform or normal distributions, more advanced methods--such as Xavier and Kaiming initialization--offer significantly improved performance, particularly in the training of deeper neural networks. Another widely used approach involves leveraging pre-trained weights from existing models to provide a more informed starting point. The key initialization strategies are discussed below.
		
		Kaiming initialization, introduced by \cite{He}, is specifically tailored for neural networks employing rectified linear unit (ReLU) activations \cite{Nair}. It was developed to address limitations of the earlier Xavier initialization \cite{Glorot}, which assumes linear activations---a condition not satisfied by ReLU-based networks. In \cite{He}, it is demonstrated that Xavier initialization could lead to vanishing or exploding gradients in deep architectures, thus impeding convergence.
		
		The central insight of Kaiming initialization is to preserve the variance of activations across layers. To achieve this, the weights of each layer are sampled from a Gaussian distribution with zero mean and a standard deviation of \( \sqrt{2/n_l} \), where \( n_l \) is the number of input connections to a neuron. This approach stabilizes the forward signal flow and accelerates convergence, particularly in very deep networks.
		
		\subsection{Neural Network Training Mechanics}
		Training a feedforward neural network involves two primary phases: the \textit{forward pass} and the \textit{backward pass}. During the forward pass, the network processes the input data to produce an output prediction. This is followed by the backward pass, where the backpropagation algorithm is employed to iteratively update the network's weights and biases based on the observed error, thereby optimizing the model.
		
		Let the training data set be represented by the matrices:
		\[
		X =
		\begin{bmatrix}
			x_1 \\
			\vdots \\
			x_m
		\end{bmatrix}
		\in \mathbb{R}^{m \times n}, \quad
		Y =
		\begin{bmatrix}
			y_1 \\
			\vdots \\
			y_m
		\end{bmatrix}
		\in \mathbb{R}^m,
		\]
		where \( m \) denotes the number of training samples and \( n \) is the input dimensionality.
		
		During network evaluation, the intermediate activations at layer \( l \) are computed and stored in the matrix:
		\[
		Z^{(l)} =
		\begin{bmatrix}
			z_1^{(l)} \\
			\vdots \\
			z_m^{(l)}
		\end{bmatrix}
		\in \mathbb{R}^{m \times n_l}.
		\]
		
		\subsection{Cost Function in Neural Network Training}
		The ultimate objective of the neural network is to learn an approximation \( \hat{f}(x) \) of the true pricing function \( f(x) \), mapping input vectors \( x \) to their corresponding outputs \( y \). To quantify the prediction error, a cost function--also referred to as the loss or error function--is utilized. This function measures the discrepancy between the network's predictions \( \hat{y} \) and the actual target values \( y \). A lower cost indicates better performance.
		
		A widely adopted loss metric is the Mean Squared Error (MSE), which computes the average of the squared differences between predicted and true values:
		\begin{equation*}
			\mathcal{C} = \frac{1}{m} \|Y - Z^{(L)}\|^2 = \frac{1}{m} (Y - Z^{(L)})^\top (Y - Z^{(L)}).
		\end{equation*}
		Here, \( Z^{(L)} \) denotes the network's output at the final layer \( L \), and the objective is to minimize \( \mathcal{C} \) through iterative optimization.

		\subsection{Gradient Descent Optimization}
		With the cost function defined, the next step in training a neural network involves minimizing this error using the gradient descent optimization technique. This iterative process adjusts the network's weights and biases to progressively improve prediction accuracy by moving in the direction that reduces the loss.
		
		To optimize the network, it is computed the gradient of the cost function with respect to each parameter. These gradients, which indicate the direction and rate of steepest ascent of the cost, are derived through partial differentiation. By moving against the gradient--i.e., in the direction of steepest descent--we ensure the cost function decreases.
		
		Consider a simplified neural network consisting of a single neuron with two weights, \( w_1 \) and \( w_2 \), and let the activation function be defined as \( \phi(y) = y^2 \). If the cost function is given as $C(z) = 2z$, where $z = w_1x_1 + w_2x_2$, then the full cost becomes:
		\begin{equation*}
			C = 2(w_1x_1 + w_2x_2)^2.
		\end{equation*}
		To find the update direction for $w_1$, we take the partial derivative:
		\begin{equation*}
			\frac{\partial C}{\partial w_1} = 4x_1(w_1x_1 + w_2x_2).
		\end{equation*}
		
		Manually computing gradients for each parameter in a deep network would be highly inefficient. Instead, backpropagation is employed--a powerful algorithm that automates and accelerates this process. Working backward from the output layer to the input layer, backpropagation applies the chain rule to efficiently propagate the gradient through each layer by reusing computations.
		
		Once the gradients are known, the parameters are updated by stepping in the direction opposite to the gradient. This update is scaled by a hyperparameter called the \textit{learning rate} ($\text{lr}$), which controls how large each adjustment is. For example, the update rule for $w_1$ becomes:
		\[
		w_1 \leftarrow w_1 - \text{lr} \cdot \frac{\partial C}{\partial w_1}
		\]
		Through successive iterations, this process steers the network toward optimal parameter values that minimize prediction error.
		
		\subsection{Stochastic Gradient Descent}
		In classical gradient descent, model training proceeds by computing the gradient of the total cost function only after the entire data set has been processed. While theoretically sound, this full-batch approach suffers from two significant drawbacks: it discards early insights from initial training samples and incurs a high computational burden, especially for large data sets.
		
		To address these limitations, Stochastic Gradient Descent (SGD) introduces a more efficient and practical alternative. Rather than relying on the complete data set, SGD approximates the true gradient by computing it over randomly selected subsets of the data, known as \textit{mini-batches}. After each mini-batch is passed through the network, the backward pass is performed and the parameters are updated accordingly.
		
		This mini-batch approach offers a trade-off between stability and stochasticity. Larger batches yield more stable and accurate gradient estimates, as they better reflect the overall data distribution. Conversely, smaller batches introduce noise into the optimization trajectory, which can be beneficial for escaping shallow local minima and improving generalization.
		
		Training progresses iteratively: each batch update constitutes one iteration. Once every sample in the data set has been used for training, a full \textit{epoch} is said to be completed. This batch-wise learning framework significantly reduces computational complexity while preserving convergence efficacy, making SGD a cornerstone of modern deep learning optimization.

		\subsection{Adam Optimizer}
		One of the earliest enhancements to the traditional gradient descent algorithm was the incorporation of \textit{momentum}, a mechanism that enables the optimization process to maintain velocity in consistent gradient directions. This not only accelerates convergence but also stabilizes parameter updates, particularly when gradients do not change direction drastically. Over time, numerous advanced optimization techniques have been proposed, among which Adam (Adaptive Moment Estimation) has become one of the most widely adopted.
		
		Adam extends the idea of momentum by simultaneously estimating both the first and second moments of the gradients. Specifically, it maintains exponentially decaying averages of past gradients $m_1$ (first moment) and their squared values $m_2$ (second moment), updated as follows:
		\begin{align*}
			m_1 &= \beta_1 \cdot m_1 + (1 - \beta_1) \cdot dx, \\
			m_2 &= \beta_2 \cdot m_2 + (1 - \beta_2) \cdot dx^2,
		\end{align*}
		where $dx$ denotes the gradient computed via backpropagation, and $\beta_1$, $\beta_2$ are hyperparameters typically set to $0.9$ and $0.999$, respectively.
		
		Since both $m_1$ and $m_2$ are initialized at zero, Adam introduces bias-corrected estimates to improve stability during the initial iterations:
		\begin{equation*}
			\hat{m}_1 = \frac{m_1}{1 - \beta_1^t}, \quad \hat{m}_2 = \frac{m_2}{1 - \beta_2^t},
		\end{equation*}
		where $t$ denotes the current iteration step. The final parameter update rule is given by:
		\begin{equation*}
			x = x - \frac{lr}{\sqrt{\hat{m}_2} + \epsilon} \cdot \hat{m}_1,
		\end{equation*}
		with $lr$ representing the learning rate and $\epsilon$ a small constant added for numerical stability.
		
		Adam's appeal lies in its efficiency, robust performance across a range of deep learning tasks, and ease of implementation. Although it requires slightly more memory to store moment estimates, this cost is negligible in most practical applications. As a result, Adam is frequently selected as the default optimizer in modern neural network training pipelines.
		
		\vspace{0.3cm}
		With the network architecture and training methodology in place, we next focus on the integration of this deep learning framework into the Heston model calibration pipeline. In the following section, we demonstrate how the trained neural networks--specifically, the PAN and CCN--are applied to enhance calibration accuracy and computational efficiency within the Heston pricing framework.
		
		Having outlined the general training strategies and optimization algorithms, we now specify the concrete training configuration adopted for PAN and CCN in the present study.
		Both PAN and CCN were trained within a supervised regression framework using MSE as the loss function and the Adam optimizer. 
		Given the limited size of the available data set, full-batch training was employed, with the batch size equal to the entire training set. Inputs and targets were standardized using statistics computed from the training data, and the same normalization parameters were applied to the validation set to ensure consistency.
		For PAN, training was conducted for 5000 epochs using Adam with a fixed learning rate of $10^{-3}$ and no learning rate scheduler.
		For CCN, training was performed for 10000 epochs using Adam with a fixed learning rate of $10^{-3}$, also without a learning rate scheduler.
		Convergence was assessed by monitoring both training and validation losses throughout optimization. Training was considered to have converged once the validation loss reached a stable plateau without exhibiting upward drift.

		\section{Implementation to the Heston Model} \label{sec:implementation}
		The calibration of the Heston stochastic volatility model is executed using a two-phase deep learning framework that efficiently leverages the regression capabilities of ANNs. This approach circumvents the computational cost of traditional pricing methods and facilitates fast, accurate estimation of model parameters from market-observed option prices.
		
		In the first phase, a synthetic data set is generated by simulating a broad and representative range of Heston model parameters--such as the initial variance \( v_0 \), long-term variance \( \theta \), mean reversion speed \( \kappa \), volatility of volatility \( \sigma \), and correlation \( \rho \)--alongside option characteristics like time to maturity and moneyness. The corresponding option prices are computed using a semi-analytical solution (e.g., the Heston closed-form formula) or a numerical pricing method, and serve as the target values for supervised learning.
		
		This data set is then split into training and validation subsets. A standard FFNN is trained on the training set to learn the nonlinear pricing function that maps model parameters and option features to option prices. Once trained, the ANN effectively approximates the pricing function, enabling rapid evaluation and gradient computation via backpropagation.
		
		In the second phase, the trained neural network is employed as a surrogate pricing model during the calibration process. The goal is to minimize the discrepancy between model-generated and market-observed option prices. This is formulated as a nonlinear least squares optimization problem over the Heston parameter space. The ANN's fast inference and differentiability ensure that the calibration is not only computationally efficient but also robust to local minima.
		
		This deep learning-based calibration strategy significantly accelerates the model fitting process and allows for real-time application in dynamic market environments.

		\subsection{Price Approximator Network}
		The first deep learning model, referred to as the PAN, is constructed to approximate the functional relationship between the strike price \( K \) and the last traded option price \( P \). Formally, the network acts as a regression function \( f: \mathbb{R} \rightarrow \mathbb{R} \), where the input is the strike price and the output is the estimated last price \( \hat{P} \).
		
		The architecture comprises multiple layers that progressively transform the input into a meaningful price prediction:
		
		\begin{enumerate}
			\item \textbf{Input Layer}: Accepts the strike price \( K \) as input.
			\item \textbf{First Hidden Layer}: Applies an affine transformation followed by the hyperbolic tangent activation:
			\[
			h^{(1)} = \tanh(W^{(1)} K + b^{(1)}),
			\]
			where \( W^{(1)} \in \mathbb{R}^{8 \times 1} \), \( b^{(1)} \in \mathbb{R}^8 \), and \( h^{(1)} \in \mathbb{R}^8 \).
			\item \textbf{Second Hidden Layer}: Applies a ReLU activation to the transformed features:
			\[
			h^{(2)} = \max(0, W^{(2)} h^{(1)} + b^{(2)}),
			\]
			where \( W^{(2)} \in \mathbb{R}^{8 \times 8} \) and \( b^{(2)} \in \mathbb{R}^8 \).
			\item \textbf{Output Layer}: Outputs the predicted price:
			\[
			\hat{P} = W^{(3)^T} h^{(2)} + b^{(3)},
			\]
			where \( W^{(3)} \in \mathbb{R}^8 \), \( b^{(3)} \in \mathbb{R} \).
		\end{enumerate}
		
		This yields the complete model:
		\[
		\hat{P} = W^{(3)^T} \max(0, W^{(2)} \tanh(W^{(1)} K + b^{(1)}) + b^{(2)}) + b^{(3)}.
		\]
		
\begin{figure}[H]
			\centering
			\includegraphics[scale=0.7]{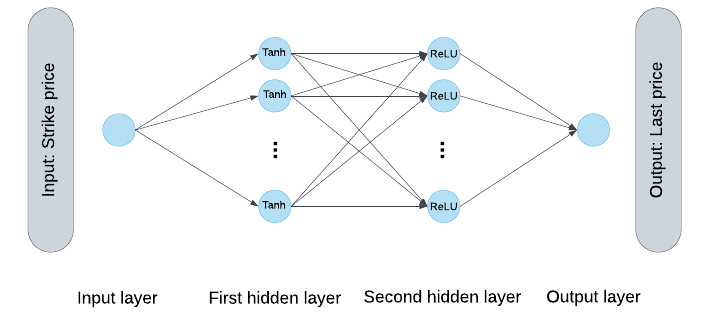}
			\caption{\label{Fig2} \small Price Approximator Network (PAN) architecture.}
		\end{figure}

		Figure~\ref{Fig2} illustrates the PAN architecture and its forward data flow.
		
		\subsection{Calibration Correction Network}
		The second model, termed the CCN, is designed to refine the pricing output of the calibrated Heston model. Given the Heston-predicted price \( P_{\text{Heston}} \), the CCN maps this value to a more accurate estimate \( \hat{P}_{\text{final}} \), effectively minimizing discrepancies between theoretical and observed market prices.
		
		Its architecture is defined as follows:
		
		\begin{enumerate}
			\item \textbf{Input Layer}: Takes \( P_{\text{Heston}} \) as input.
			\item \textbf{First Hidden Layer}: Applies a sigmoid activation:
			\[
			h^{(1)} = \sigma(W^{(1)} P_{\text{Heston}} + b^{(1)}),
			\]
			with \( W^{(1)} \in \mathbb{R}^{7 \times 1} \), \( b^{(1)} \in \mathbb{R}^7 \).
			\item \textbf{Second Hidden Layer}: Applies the \( \tanh \) activation:
			\[
			h^{(2)} = \tanh(W^{(2)} h^{(1)} + b^{(2)}),
			\]
			with \( W^{(2)} \in \mathbb{R}^{7 \times 7} \), \( b^{(2)} \in \mathbb{R}^7 \).
			\item \textbf{Output Layer}: Produces the corrected price estimate:
			\[
			\hat{P}_{\text{final}} = W^{(3)^T} h^{(2)} + b^{(3)},
			\]
			where \( W^{(3)} \in \mathbb{R}^7 \), \( b^{(3)} \in \mathbb{R} \).
		\end{enumerate}
		
		The complete expression for the CCN becomes:
		\[
		\hat{P}_{\text{final}} = W^{(3)^T} \tanh(W^{(2)} \sigma(W^{(1)} P_{\text{Heston}} + b^{(1)}) + b^{(2)}) + b^{(3)}.
		\]
		
\begin{figure}[H]
			\centering
			\includegraphics[scale=0.17]{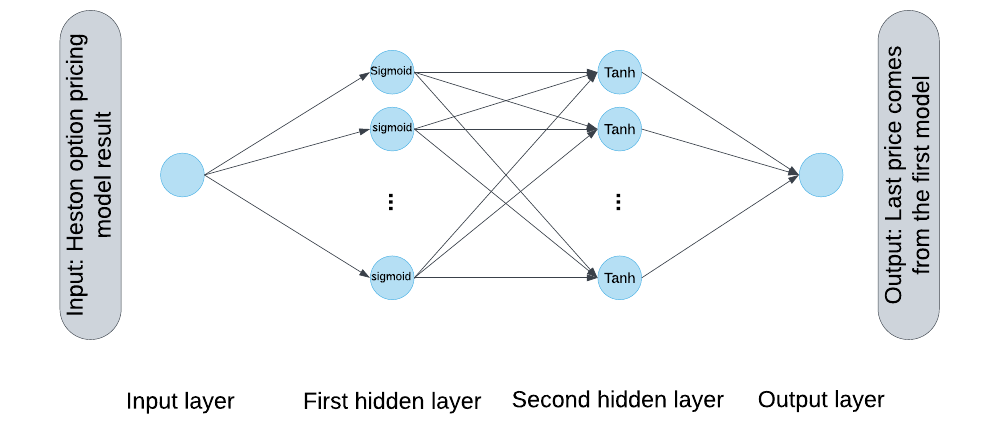}
			\caption{\label{Fig3} \small Calibration Correction Network (CCN) architecture.}
		\end{figure}

		Figure~\ref{Fig3} shows the architecture and data flow within the CCN.
		
		Traditional calibration methods often fail to accurately capture the nonlinear relationships between parameters in the Heston model, especially under volatile market regimes. Deep learning models such as PAN and CCN offer flexible function approximation capabilities and adapt to complex data structures without requiring explicit assumptions.
		
		The PAN provides a reliable proxy for market data, while the CCN enhances theoretical model alignment. Together, these networks form a robust and efficient calibration pipeline, improving parameter estimation and ultimately enabling more accurate option pricing in practice.
		
		\subsection{Activation Function Selection}
		
		The activation functions employed in the PAN and CCN were selected in light of the distinct functional roles of these two components within the calibration framework. Although both networks are shallow fully connected architectures, they serve fundamentally different objectives: PAN approximates the option pricing surface, whereas CCN refines the baseline Heston-implied prices through systematic correction.
		
		For PAN, which learns the mapping between strike price and observed option values, the target function is smooth, convex, and exhibits varying curvature across moneyness regions. To capture these characteristics effectively, a mixed activation structure was adopted. The first hidden layer uses the $\tanh$ activation, which is symmetric and continuously differentiable, making it well suited for modeling smooth global nonlinearities. Because inputs are standardized prior to training, $\tanh$ operates in a numerically stable regime and captures the overall curvature of the pricing surface efficiently. The second hidden layer employs the ReLU activation, introducing piecewise-linear flexibility that enables the representation of localized slope variations, particularly in the at-the-money region where price sensitivity changes more rapidly. In addition, ReLU facilitates stable gradient propagation and mitigates potential vanishing-gradient effects associated with purely saturating activations. The $\tanh$--ReLU combination therefore balances smooth global approximation with localized expressiveness.
		
		In contrast, CCN is designed to learn a correction function applied to Heston-implied prices. Rather than reconstructing the full pricing surface, it produces controlled residual adjustments. Consequently, bounded and smooth activation functions are more appropriate. The first hidden layer employs a sigmoid activation to constrain intermediate representations within a finite range, thereby promoting numerical stability and preventing extreme correction responses. The second hidden layer uses $\tanh$, allowing symmetric positive and negative adjustments while maintaining smoothness. Since the correction may increase or decrease the baseline model price, a symmetric activation is appropriate at this stage. The sigmoid--$\tanh$ configuration thus ensures that the correction remains stable, smooth, and economically coherent.
		
		Preliminary experiments with uniform activation schemes (e.g., exclusively ReLU or exclusively $\tanh$) did not yield statistically significant improvements in RMSE or MAE and, in some cases, resulted in slightly less stable training dynamics. Accordingly, the mixed activation configuration was retained to achieve an appropriate balance between expressive capacity, numerical stability, and controlled correction behavior.
		
		\subsection{Mathematical Framework for Deep Learning-Based Calibration}
		To formalize the calibration procedure, we define the objective function as the weighted mean squared error between model-predicted and market-observed option prices. Given a set of market data indexed by maturity \( T_i \) and log-moneyness \( m_j \), and a neural network with trained weights \( \hat{W} \), the calibration problem can be expressed as:
		\begin{equation*}
			\hat{\eta} = \arg \min_{\eta \in \mathcal{H}} \sum_{i,j} \omega_{i,j} \left( F(\hat{W}, \eta_i, T_i, m_j) - V^{\text{Mkt}}(T_i, m_j) \right)^2.
		\end{equation*}
		The proposed framework leverages two neural networks to improve pricing accuracy. The PAN is trained to approximate the mapping between strike price \( K \) and the last traded option price \( \hat{P}_{\text{PAN}}(K) \). Subsequently, the CCN refines the output of the calibrated Heston model by learning a correction function:
		
		\begin{equation*}
			\hat{P}_{\text{final}} = \text{CCN}(P_{\text{Heston}}).
		\end{equation*}
		
		This two-stage architecture enables the model to first emulate the option price surface from data (via PAN) and then enhance calibration precision by correcting discrepancies (via CCN), even when traditional calibration is suboptimal.
		
		To quantify the performance of both approaches, we compute standard error metrics including Root Mean Squared Error (RMSE), Mean Absolute Error (MAE), and Mean Relative Error (MRE), defined respectively as:
		
		\begin{align*}
			\text{RMSE} &= \sqrt{ \frac{1}{N} \sum_{i=1}^{N} \left( P^{\text{model}}_i - P^{\text{market}}_i \right)^2 }, \\
			\text{MAE} &= \frac{1}{N} \sum_{i=1}^{N} \left| P^{\text{model}}_i - P^{\text{market}}_i \right|, \\
			\text{MRE} &= \frac{1}{N} \sum_{i=1}^{N} \frac{\left| P^{\text{model}}_i - P^{\text{market}}_i \right|}{P^{\text{market}}_i}.
		\end{align*}
		
		These metrics provide a comprehensive assessment of model accuracy across different pricing levels and moneyness intervals. The following figures and tables illustrate how the proposed deep learning framework achieves lower calibration errors compared to the traditional method.

		\section{Empirical Results and Analysis} \label{results}
		This section presents a comprehensive evaluation of the proposed two-phase deep learning framework for calibrating the Heston stochastic volatility model. The empirical results are analyzed to assess both the effectiveness and robustness of the methodology, with a focus on comparing its performance to that of traditional calibration techniques.
		
		The evaluation is conducted using both in-sample and out-of-sample data sets to ensure the generalizability of the model. We systematically examine the calibration accuracy by comparing model-implied prices with observed market prices. To quantify performance, we employ standard error metrics such as MSE and MAE, which provide insight into the predictive accuracy and consistency of the deep learning-enhanced Heston model.
		
\begin{figure}[H]
			\centering
			{\includegraphics[scale=0.6]{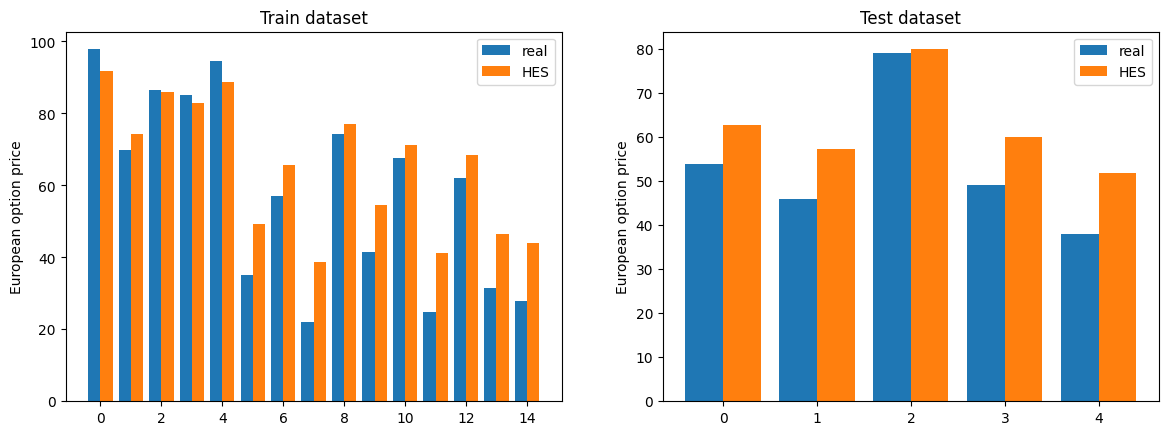}}
			\caption{\label{Fig4} \small Results of the calibration of the Heston's option pricing model on $S\&P$ 500 option data.}
		\end{figure}

		As part of the empirical study, we analyze European option prices on the S\&P 500 index as of February 7, 2025, with the underlying asset priced at \$6025.99. The analysis focuses on near-term options expiring within three days, utilizing the latest available market data from the final trading session. The calibration results are presented in Figure~\ref{Fig4}, which highlights the observed deviations between market option prices and those produced by the Heston model under conventional calibration methods.
		
		These discrepancies underscore the limitations of traditional approaches in capturing the complex dynamics of market prices, particularly under short-term volatility conditions. The findings demonstrate the necessity for advanced calibration strategies--such as those employing deep neural networks--to reduce pricing errors and enhance the model's alignment with real-world market behavior.
		
\begin{figure}[H]
			\centering
			{\includegraphics[scale=0.55]{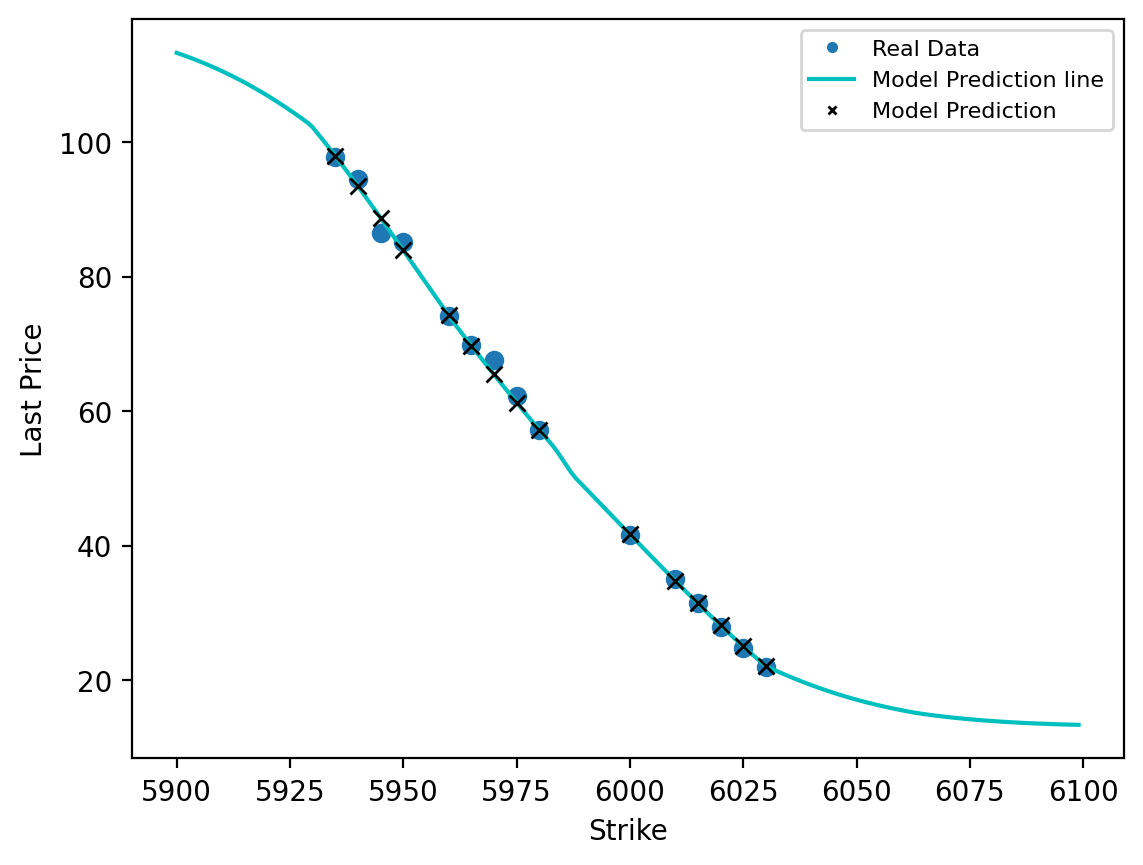}}
			\caption{\label{Fig5} \small Approximation of a line using the Price Approximator Network (PAN) model on $S\&P$ 500 option.}
		\end{figure}

\begin{figure}[H]
			\centering
			{\includegraphics[scale=0.6]{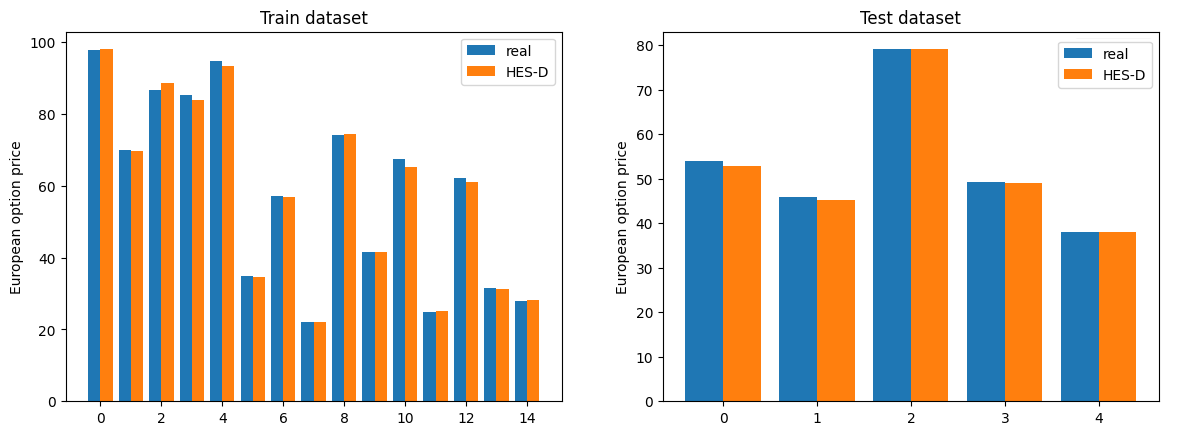}}
			\caption{\label{Fig6} \small Results of improving the Heston's option pricing model calibration using deep learning for $S\&P$ 500 data.}
		\end{figure}

		The PAN is first employed to estimate a smooth pricing curve for the last traded prices of $S\&P$ 500 options, using in-sample data. The results of this estimation are depicted in Fig. \ref{Fig5}, which demonstrates the network's ability to accurately capture the underlying pricing structure. Building on this approximation, the CCN is subsequently trained using both the PAN-derived price surface and the preliminary calibration outputs from the Heston model. The goal is to refine these outputs and reduce systematic discrepancies. The enhanced pricing performance achieved through this correction mechanism is illustrated in Fig. \ref{Fig6}, clearly indicating that the proposed deep learning-enhanced calibration framework outperforms the traditional approach. To further substantiate this improvement, error metrics such as MSE and MAE are reported in Table~\ref{tab:calibration_results}, validating the superior performance of the deep learning--augmented approach over conventional calibration techniques.
		
		%
		
		\begin{table}[H]
			\centering
			\caption{Error metrics for $S\&P$ 500 data calibration.}
			\label{tab:calibration_results}
			\begin{tabular}{lrr}
				\toprule
				\textbf{Metric} & \textbf{Traditional} & \textbf{Deep Learning} \\
				\midrule
				Train RMSE      & 10.47                & 0.98                   \\
				Train MRE       & 0.2515               & 0.0101                 \\
				Train MAE       & 8.81                 & 0.65                   \\
				Test RMSE       & 10.21                & 0.55                   \\
				Test MRE        & 0.2023               & 0.0079                 \\
				Test MAE        & 9.20                 & 0.39                   \\
				\bottomrule
			\end{tabular}
		\end{table}
		
		As demonstrated in Table \ref{tab:calibration_results}, the deep learning-based framework significantly outperforms the traditional calibration method across all error metrics. On the training set, the deep learning approach achieves notably lower values for RMSE ($0.98$ vs.\ $10.47$), MRE ($0.0101$ vs.\ $0.2515$), and MAE ($0.65$ vs.\ $8.81$), indicating a much closer alignment with observed market prices. The model maintains this superior performance on the test set, achieving an RMSE of $0.55$, substantially lower than the traditional method's $10.21$, alongside similar improvements in MRE and MAE. These results affirm the effectiveness of the proposed approach in refining the calibration of the Heston option pricing model.
		
\begin{figure}[H]
			\centering
			{\includegraphics[scale=0.6]{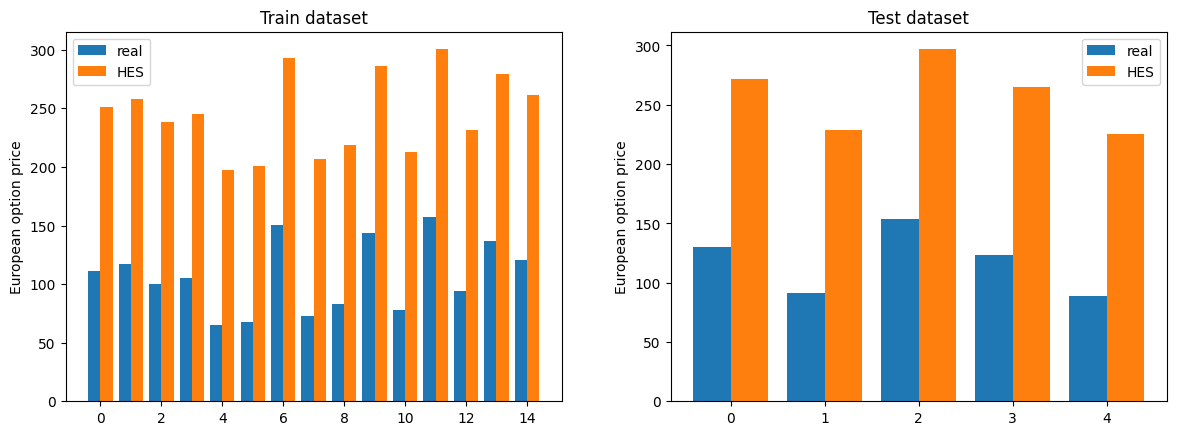}}
			\caption{\label{Fig7} \small Results of the calibration of the Heston's option pricing model on $S\&P$ 500 Mini option data.}
		\end{figure}

		To assess the generalizability of the deep learning framework, we extend the analysis to a different data set--$S\&P$ 500 Mini options with 39 days to expiration, evaluated as of March 20, 2025 (with calibration conducted on February 7, 2025). At this point, the underlying index value is \$6049.50. We first apply the traditional calibration approach to this new data set, with the resulting model-generated option prices presented in Fig.~\ref{Fig7}. This step provides a baseline for comparison, particularly in scenarios involving longer maturities and different market regimes.
		
\begin{figure}[H]
			\centering
			{\includegraphics[scale=0.55]{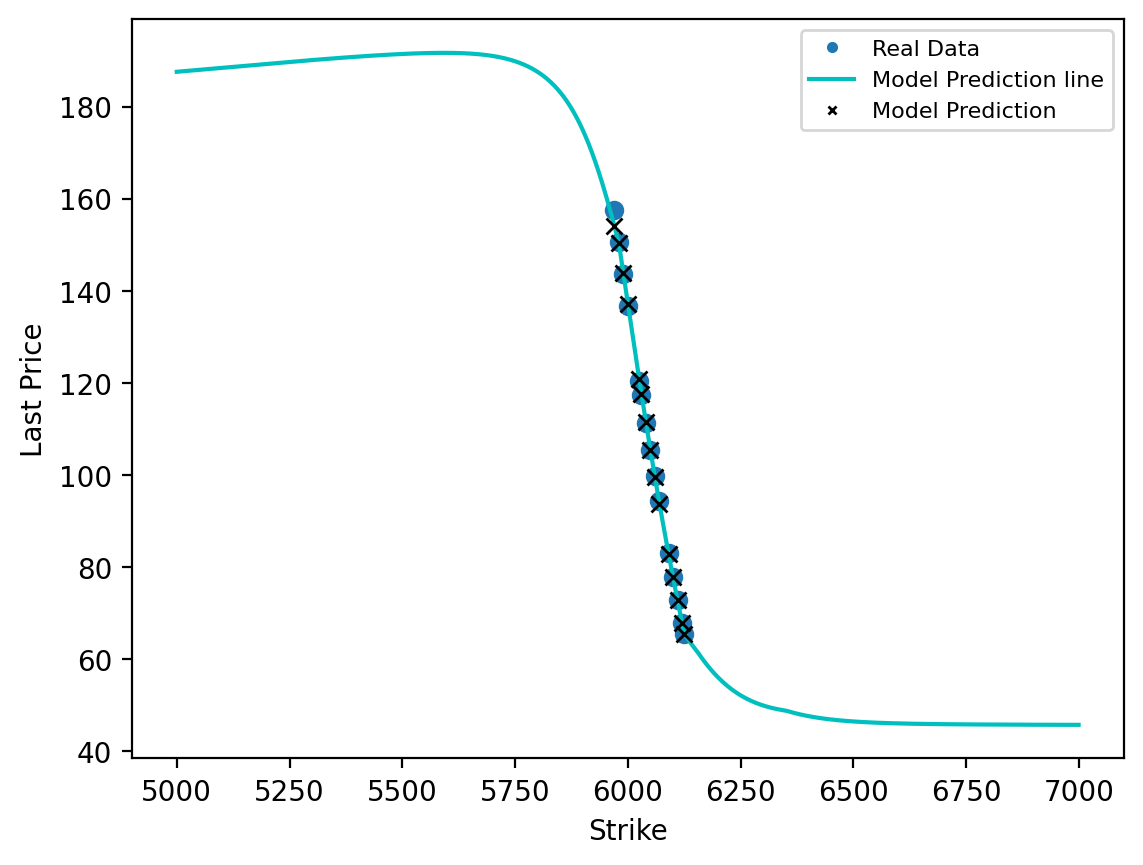}}
			\caption{\label{Fig8} \small Approximation of a line using the Price Approximator Network (PAN) model on $S\&P$ 500 mini option data.}
		\end{figure}

\begin{figure}[H]
			\centering
			{\includegraphics[scale=0.6]{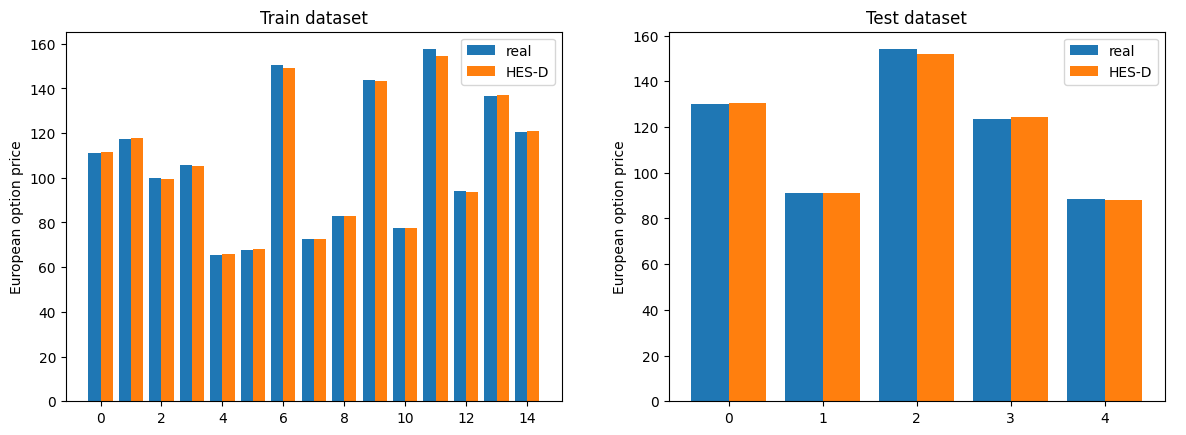}}
			\caption{\label{Fig9} \small Results of improving the Heston's option pricing model calibration using deep learning for $S\&P$ 500 mini data.}
		\end{figure}

		Next, the PAN is applied to the $S\&P$ 500 Mini data set to learn a smoothed representation of the observed option price surface (see Fig.~\ref{Fig8}). Building upon this, the CCN is trained to adjust the outputs of the initial Heston calibration, resulting in a refined price estimation. The improved performance of this combined method is illustrated in Fig.~\ref{Fig9}. For a quantitative evaluation, Table \ref{tab:calibration_results1} compares the deep learning-enhanced approach with the traditional method across multiple accuracy metrics. The results clearly show that the deep learning model maintains its superiority. On the training set, the RMSE, MRE, and MAE are dramatically reduced to $0.89$, $0.005$, and $0.56$ respectively, compared to the traditional method's $138.60$, $1.39$, and $138.55$. This performance gap persists in the test set, where the deep learning model achieves an RMSE of $1.02$ versus $140.11$ for the traditional method. These findings underscore the robustness and adaptability of the proposed deep learning framework, particularly in capturing option price dynamics across varying maturities and market environments.
		
		%
		
		\begin{table}[H]
			\centering
			\caption{Error metrics for $S\&P$ 500 Mini data calibration.}
			\label{tab:calibration_results1}
			\begin{tabular}{lrr}
				\toprule
				\textbf{Metric} & \textbf{Traditional} & \textbf{Deep Learning} \\
				\midrule
				Train RMSE      & 138.60               & 0.89                   \\
				Train MRE       & 1.39                 & 0.005                  \\
				Train MAE       & 138.55               & 0.56                   \\
				Test RMSE       & 140.11               & 1.02                   \\
				Test MRE        & 1.24                 & 0.006                  \\
				Test MAE        & 140.09               & 0.78                   \\
				\bottomrule
			\end{tabular}
		\end{table}
		
			
			\subsection{Extended Out-of-Sample Validation and Error Analysis}
			
			To further evaluate the robustness of the proposed PAN--CCN framework, we conducted additional out-of-sample experiments on assets beyond the S\&P 500. Specifically, we evaluated the methodology on NASDAQ-100 index options and VIX options, which represent different underlying structures and pricing scales.
			
			We consider European call options on the NASDAQ-100 index as of 21 February 2026, with expiration on 2 March 2026 (10 days to maturity). The underlying index level at the evaluation date was 24,984.65. The quantitative results are reported in Table~\ref{tab:nasdaq_results}.
			
			\begin{table}[H]
				\centering
				\caption{Out-of-sample calibration performance on NASDAQ-100 European call options (21 February 2026; maturity: 2 March 2026).}
				\label{tab:nasdaq_results}
				\begin{tabular}{lcccccc}
					\hline
					& \multicolumn{3}{c}{Train} & \multicolumn{3}{c}{Test} \\
					\cline{2-4} \cline{5-7}
					Method & RMSE & MAE & MRE & RMSE & MAE & MRE \\
					\hline
					Classical Heston 
					& 149.8603 & 147.7591 & 2.1805 
					& 126.4693 & 117.8126 & 4.3027 \\
					
					PAN--CCN 
					& 9.7467 & 5.6477 & 0.0470 
					& 48.3769 & 27.1886 & 0.2820 \\
					\hline
				\end{tabular}
			\end{table}
			
			Under classical Heston calibration, the test RMSE is 126.4693 and the test MAE is 117.8126. After applying the PAN--CCN correction, the test RMSE decreases to 48.3769 and the test MAE to 27.1886. This corresponds to a substantial reduction in absolute pricing errors. A similar improvement is observed in relative error (MRE), which decreases from 4.3027 to 0.2820 on the test set.
			
\begin{figure}[H]
				\centering
				\includegraphics[width=0.48\textwidth]{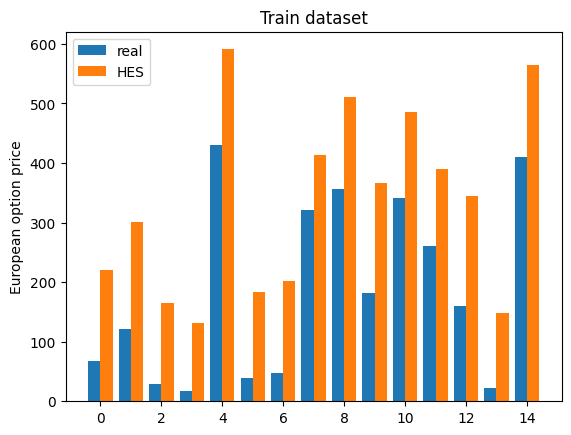}
				\includegraphics[width=0.48\textwidth]{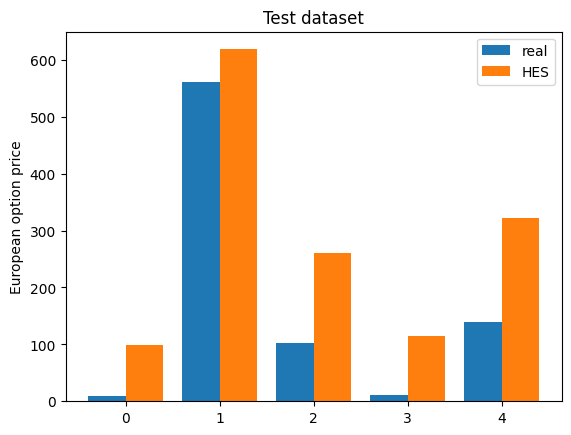}
				\caption{Baseline Heston calibration results on NASDAQ-100 European option data (21 Feb 2026; maturity: 2 Mar 2026). The plots compare market prices with Heston-implied prices.}
				\label{fig:nasdaq_baseline}
			\end{figure}

			Figure~\ref{fig:nasdaq_baseline} shows that, under classical calibration, noticeable deviations between model and market prices occur across several strikes in the test sample. In particular, the largest discrepancies appear in higher-priced contracts, where the absolute deviations are visually larger. After applying the neural network correction (Figure~\ref{fig:nasdaq_dl}), the gap between model and market prices narrows visibly across all reported strikes. While residual errors remain, their magnitude is consistently smaller than in the baseline case, which is consistent with the quantitative improvements reported in Table~\ref{tab:nasdaq_results}.
			
			Since the reported strikes span in-the-money, at-the-money, and out-of-the-money regions, the graphical comparisons provide a visual assessment of error behavior across different moneyness levels. The reduction in deviations after applying PAN--CCN is observed across the strike range rather than being confined to a narrow subset of contracts.
			
\begin{figure}[H]
				\centering
				\includegraphics[width=0.75\textwidth]{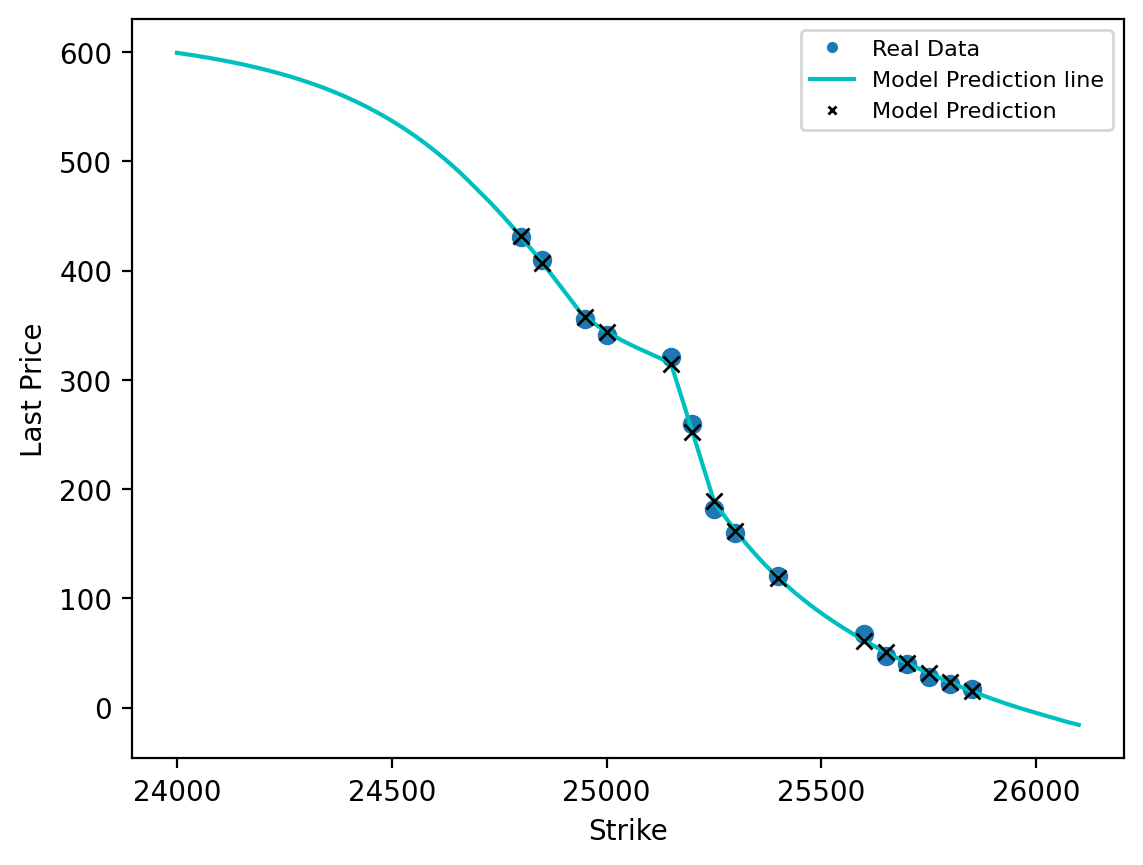}
				\caption{PAN approximation of the NASDAQ-100 option price curve as a function of strike. Circles denote market prices, while the curve represents the PAN approximation.}
				\label{fig:nasdaq_pan}
			\end{figure}

\begin{figure}[H]
				\centering
				\includegraphics[width=0.48\textwidth]{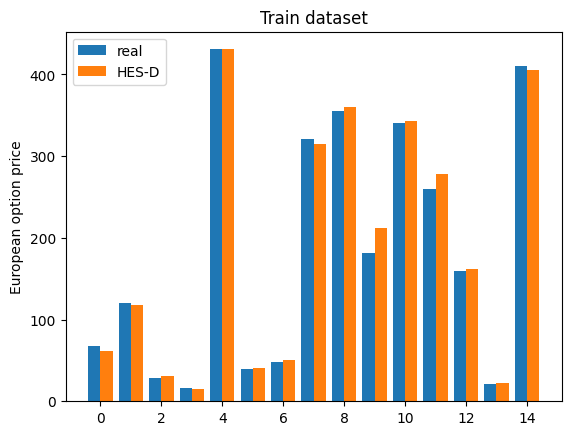}
				\includegraphics[width=0.48\textwidth]{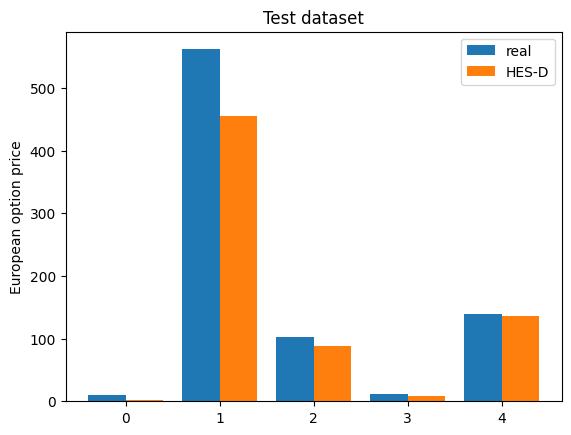}
				\caption{Deep learning--enhanced calibration on NASDAQ-100 option data. The CCN-corrected prices show improved agreement with market prices relative to the baseline calibration.}
				\label{fig:nasdaq_dl}
			\end{figure}

			The PAN approximation shown in Figure~\ref{fig:nasdaq_pan} exhibits a smooth and monotonic relationship between strike and option price. The corrected pricing curve remains visually smooth and does not display irregular local fluctuations.
			
			We additionally consider VIX European call options as of 21 February 2026, with expiration on 15 April 2026 (53 days to maturity). The VIX level at the evaluation date was 19.09. The corresponding results are summarized in Table~\ref{tab:vix_results}.
			
			\begin{table}[H]
				\centering
				\caption{Out-of-sample calibration performance on VIX European call options (21 February 2026; maturity: 15 April 2026).}
				\label{tab:vix_results}
				\begin{tabular}{lcccccc}
					\hline
					& \multicolumn{3}{c}{Train} & \multicolumn{3}{c}{Test} \\
					\cline{2-4} \cline{5-7}
					Method & RMSE & MAE & MRE & RMSE & MAE & MRE \\
					\hline
					Classical Heston 
					& 0.2479 & 0.1707 & 0.0304 
					& 0.2103 & 0.1736 & 0.0389 \\
					
					PAN--CCN 
					& 0.1986 & 0.1252 & 0.0262 
					& 0.1447 & 0.0907 & 0.0159 \\
					\hline
				\end{tabular}
			\end{table}
			
			For classical calibration, the test errors are RMSE = 0.2103, MAE = 0.1736, and MRE = 0.0389. Using the PAN--CCN framework, the test errors improve to RMSE = 0.1447, MAE = 0.0907, and MRE = 0.0159 (Table~\ref{tab:vix_results}).
			
\begin{figure}[H]
				\centering
				\includegraphics[width=0.48\textwidth]{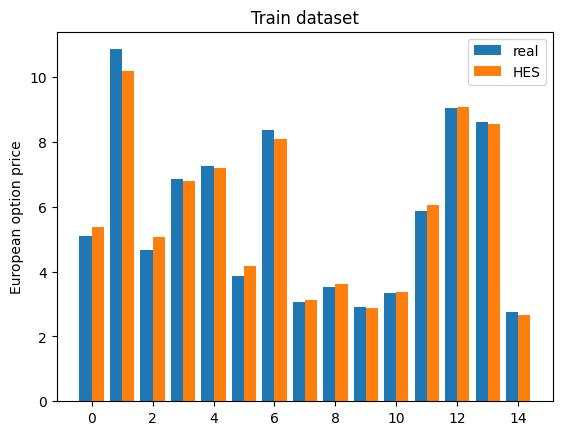}
				\includegraphics[width=0.48\textwidth]{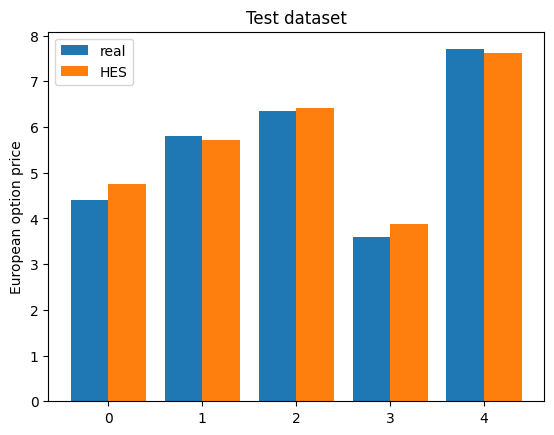}
				\caption{Baseline Heston calibration results on VIX European option data (21 Feb 2026; maturity: 15 Apr 2026). The plots compare market prices with Heston-implied prices.}
				\label{fig:vix_baseline}
			\end{figure}

\begin{figure}[H]
				\centering
				\includegraphics[width=0.75\textwidth]{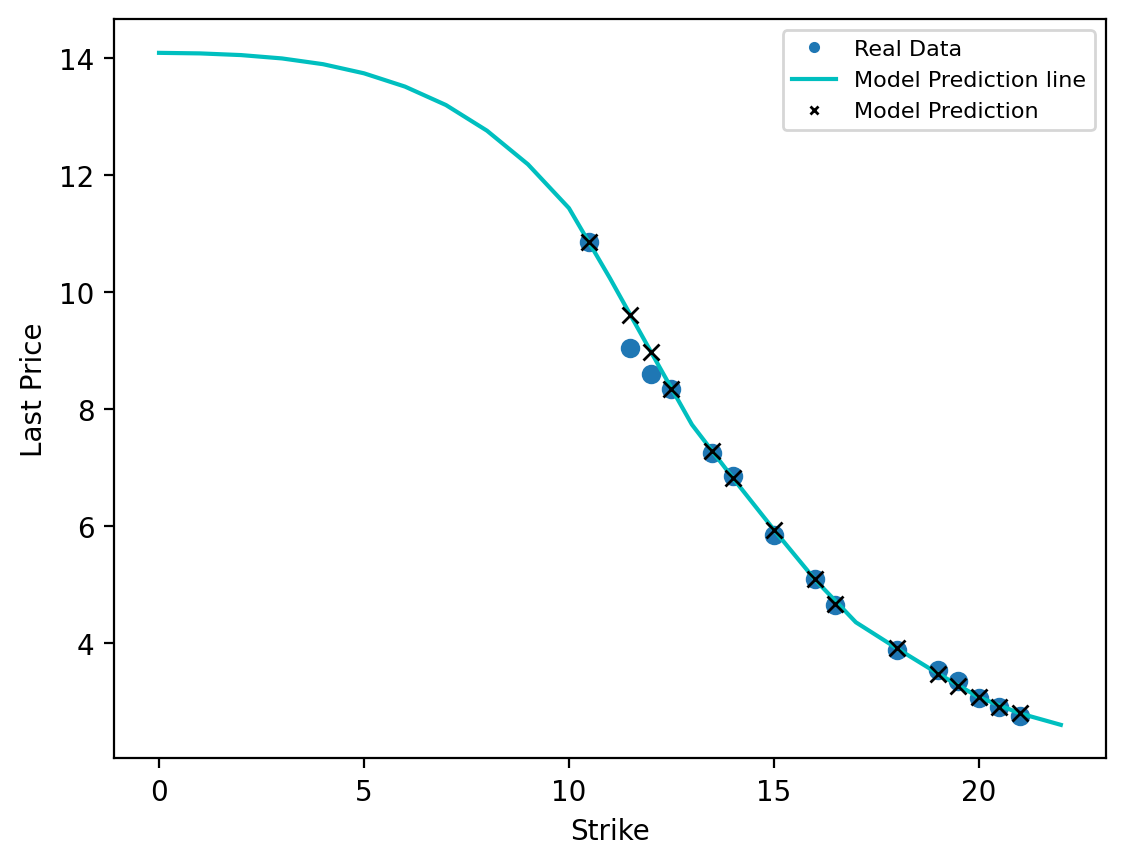}
				\caption{PAN approximation of the VIX option price curve as a function of strike. Circles denote market prices, while the curve represents the PAN approximation.}
				\label{fig:vix_pan}
			\end{figure}

\begin{figure}[H]
				\centering
				\includegraphics[width=0.48\textwidth]{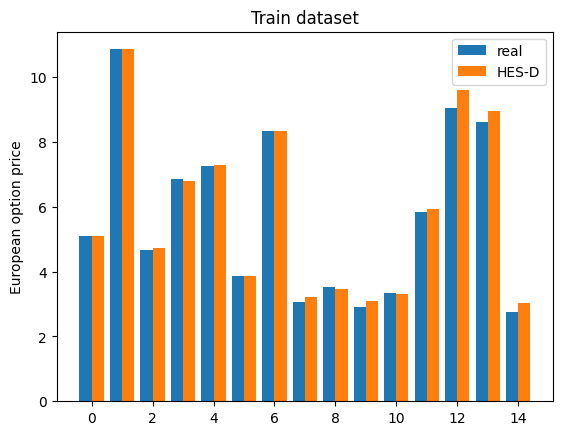}
				\includegraphics[width=0.48\textwidth]{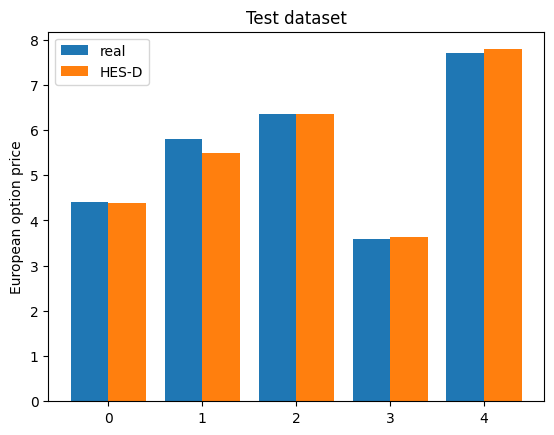}
				\caption{Deep learning--enhanced calibration on VIX option data. The CCN-corrected prices exhibit improved alignment with market prices compared to the baseline calibration.}
				\label{fig:vix_dl}
			\end{figure}

			Figure~\ref{fig:vix_baseline} compares market prices and Heston-implied prices for the VIX sample, and Figure~\ref{fig:vix_dl} shows the corresponding comparison after deep learning enhancement. In line with the reductions reported in Table~\ref{tab:vix_results}, the corrected outputs in Figure~\ref{fig:vix_dl} are visually closer to the market prices than the baseline results. Figure~\ref{fig:vix_pan} presents the PAN approximation as a smooth fitted curve over strikes for this sample.
			
			Although each asset is evaluated at a fixed maturity, the two datasets correspond to materially different time-to-maturity settings (10 days for NASDAQ-100 and 53 days for VIX), providing validation under distinct horizon lengths. Across both assets and their respective strike ranges, the improvements are not localized to a single pricing region but are observed broadly across the reported contracts.
			
			Taken together, the numerical evidence in Tables~\ref{tab:nasdaq_results}--\ref{tab:vix_results} and the graphical comparisons in Figures~\ref{fig:nasdaq_baseline}--\ref{fig:vix_dl} support the conclusion that the proposed framework improves out-of-sample calibration accuracy across distinct underlying assets and pricing regimes, while maintaining the overall smooth shape of the option price curve as illustrated in the figures.
			
		\section{Conclusion} \label{conclusion}
		In this study, we proposed a novel deep learning-based framework to enhance the calibration accuracy of the Heston stochastic volatility model. The methodology incorporated two regression-based neural networks: the Price Approximator Network (PAN), which approximates the relationship between strike prices and option values, and the Calibration Correction Network (CCN), which systematically refines baseline Heston-implied prices to better align with market observations.
		
		The primary objective of the proposed calibration procedure was to achieve superior out-of-sample performance while preserving computational efficiency. Extensive empirical analyses conducted across heterogeneous assets—including equity index options (NASDAQ-100) and volatility index derivatives (VIX)—demonstrated that the deep learning-enhanced framework consistently reduced pricing errors relative to classical calibration. The improvements were observed in both absolute (RMSE, MAE) and relative (MRE) error metrics and remained stable across distinct market regimes.
		
		Importantly, error reduction was not confined to isolated regions of the option surface. The proposed architecture exhibited improved performance across different moneyness levels and maturities, including high-curvature at-the-money regions and short-maturity contracts, where traditional calibration methods often struggle. This confirms that the learned correction function generalizes smoothly rather than overfitting localized strike regions.
		
		By effectively capturing complex nonlinear dependencies in the option pricing space, the proposed PAN--CCN framework provides a robust and computationally efficient alternative to conventional optimization-based calibration techniques. The cross-asset validation results further suggest that the architecture is model-agnostic in structure and potentially extendable to other stochastic volatility and jump-diffusion frameworks.
		
		Overall, the findings indicate that deep learning can serve not merely as a surface approximation tool, but as a structured enhancement layer that systematically improves model-based calibration across diverse financial environments.

%
%


	\end{document}